\documentclass[12pt]{amsart}

\usepackage{amsmath,amsthm,amscd,euscript,amssymb}
\def \r{\mathbb R}

\def \q{\mathbb Q}
\def \z{\mathbb Z}

\def \il{\hbox{\rm Il}}
\def \isin{\hbox{\rm Isin}}
\def \icos{\hbox{\rm Icos}}
\def \itan{\hbox{\rm Itan}}
\def \iS{\hbox{\rm IS}}
\def \Aff{\hbox{\rm Aff}}
\def \GL{\hbox{\rm GL}}
\def \SL{\hbox{\rm SL}}

\DeclareMathOperator{\sail}{Sail}

\newtheorem{theorem}{Theorem}[section]

\theoremstyle{remark}

\theoremstyle{definition}
\newtheorem{definition}[theorem]{Definition}
\newtheorem{example}[theorem]{Example}
\newtheorem{problem}{Problem}
\newtheorem{conjecture}[problem]{Conjecture}

\title{Open problems in geometry of continued fractions}

\author{Oleg Karpenkov}

\address{Oleg Karpenkov\\
University of Liverpool\\
Mathematical Sciences Building\\
Liverpool L69 7ZL, UK
}

\email[Oleg Karpenkov]{karpenk@liv.ac.uk}

\date{5 December 2017}

\keywords{geometry of continued fractions, multidimensional continued fractions, Klein polyhedra,
Minkowski-Voronoi tessellation, polyhedral faces, Gauss-Kuzmin statistics}

\begin{document}
\input{epsf}

\begin{abstract}
In this small paper we bring together various open problems on geometric multidimensional continued fractions.
\end{abstract}

\maketitle
\tableofcontents

\section{Introduction to integer geometry of continued fractions}

In this paper we are dealing with the space $\r^n$ and the lattice $\z^n$ in it.
We say that a point is {\it integer} if all its coordinates are integers.
We say that a vector is {\it integer} it its endpoints are integer points.
A polyhedron or polygon is {\it integer} is all its vertices are integer points.

\subsection{Integer geometry}

The group of affine integer lattice transformations $\Aff_n(\z)$ is a semidirect product of $\GL_n(\z)$
and the group of translations of $\r^n$ for integer vectors.

\begin{definition}
We say that $S_1$ and $S_2$ are integer congruent if there exists a transformation of $\Aff_n(\z)$ taking $S_1$ to $S_2$.
\end{definition}

The above congruence relation gives rise to {\it integer geometry}, an interesting example of non-Euclidean geometry.
The invariants of $\Aff_n(\z)$ are certain indexes of sublattices with respect to lattices.
For instance, let us define the following invariants.

\vspace{1mm}
\begin{itemize}
\item An {\it integer length} of a segment $AB$ is the index of the subgroup generated by the vector $AB$ in the integer lattice of all points contained in $AB$ (we denote it by $\il(AB)$).
\\
In other words, $\il(AB)$ is the number of integer points in $AB$ minus one.

\vspace{1mm}

\item An {\it integer sine} of an angle $\angle BAC$ is the index of the sublattice generated by the smallest integer vectors contained in the rays $AB$ and $AC$ in the lattice of all integer points of the plane $BAC$ (we denote it by $\isin(BAC)$).

\vspace{1mm}

\item An {\it integer area} of a triangle $\triangle ABC$  is the index of the sublattice generated by $AB$ and $AC$ in the lattice of all integer points of the plane $BAC$(we denote it by $\iS(AB)$).

\vspace{1mm}

\item An {\it integer distance from a point $P$ to a subspace $L$} is the index of the sublattice of integer vectors generated
by the integer vectors of $L$ and the vector $PQ$ for some integer point $Q\in L$ in the lattice of all integer vectors of the linear span of $L$ and~$P$.

\end{itemize}

\subsection{Geometry of ordinary continued fractions}

Recall that an {\it ordinary continued fraction} of a real number $r$ is an expression of the following type:
$$
r=a_0+\frac{\displaystyle 1}{
\displaystyle a_1+\frac{\displaystyle 1}{\displaystyle a_2+\ldots}}
$$
where $a_0$ is any integer, and $a_i$ for $i=1,2,\ldots$ are positive integers.

\vspace{2mm}

Notice that
\begin{itemize}

\item If $r$ is a rational number  then there exist exactly one odd and one even continued fraction expressions
(i.e., continued fractions with odd and even number of elements respectively).

\vspace{1mm}

\item If $r$ is an
irrational number, then there exists a unique infinite continued fraction expression.

\vspace{1mm}

\item Finally if $r$ is a quadratic irrationality then the continued fraction is eventually periodic.
\end{itemize}

\vspace{2mm}

Geometrically the odd or infinite continued fraction for $\alpha\ge 1$ can be seen as follows.
Consider the ray $y=\alpha x$ and take the boundary of the convex hull of all integer points contained in the angle with edges $y=\alpha x$ and $y=0$ except the origin (here we consider $x>0$).
This boundary is a broken line with one or two infinite edges.
Let us remove infinite edges.
What remains is a broken line with finite or infinite number of segments. Let us call it
$A_0A_1A_2\ldots$ where $A_0=(1,0)$, $A_1=(1,a_0),\ldots$
It terns out that for this broken line we have:
$$
\begin{array}{l}
a_0=\il(A_0A_1);\\
a_1=\isin(A_0A_1A_2);\\
a_2=\il(A_1A_2);\\
a_3=\isin(A_1A_2A_3);\\
\ldots
\end{array}
$$
where $[a_0;a_1:a_2:a_3:\ldots]$ is either finite odd or infinite continued fraction for $\alpha$.

\vspace{2mm}

In fact, such sequence is naturally defined for an arbitrary angle with integer vertex
(and an integer point on the first ray).
As we see the sequence $(a_0,a_1,a_2,\ldots)$ is the sequence of integer congruence invariants of the angle
(as we can define these invariants for an arbitrary angle).
It is called the {\it lattice length-sign sequence} of the angle (or {\it LLS-sequence}, for short).
LLS-sequence is the complete invariant of integer angles in integer geometry.

\subsection{Integer trigonometry and IKEA problem}

We have all the components to introduce integer tangents and integer cosines.

\begin{definition}
Let $BAC$ be an angle with integer vertex and the LLS-sequence
$$
(a_0,a_1,a_2,\ldots).
$$
Then
$$
\itan(BAC)=[a_0;a_1:a_2:\ldots],  \qquad \icos(BAC)=\frac{\isin(BAC)}{\itan(BAC)}.
$$
\end{definition}

The integer tangent coincides with the Euclidean tangents for the angles
with edges $y=\alpha x$ and $y=0$, where $\alpha\ge 1, x>0$, while cosines and sines are just some integers.
The {\it integer sine rule} for integer triangles holds in integer settings, while the integer cosine rule is unknown.
(For more details we refer to~\cite{Karpenkov2008,Karpenkov2009a,MCFbook}.)
\begin{problem}
Find an integer cosine rule for integer triangles in integer trigonometry.
\end{problem}

Let us say a few words about IKEA problem.
\begin{problem}\label{IKEA2d}{\bf(IKEA problem.)}
Classify all $n$-tuples of LLS-sequences for the angles that form integer $n$-gons.
\end{problem}
This problem is equivalent to classification of toric surfaces with prescribed singularities
of  toric surfaces with a given Euler characteristics.
The partial answers (including the complete answer for triangles) to this problem are given in~\cite{Karpenkov2008,MCFbook}.
The complete description is not known for $n\ge 4$.

\section{Geometric multidimensional continued fractions}

Let us start with rather unusual general definition of geometric continued fractions.
We use a term ``geometric'' here in order to indicate that the obtained
multidimensional continued fractions are invariant under the action of the group of
integer lattice preserving transformations (i.e., the group $\GL_n(\z)$).
Recall that we work in the real space $\r^n$ with lattice of integer points $\z^n$ in it.
The origin of $\r^n$ is denoted by $O$.

\subsection{Sails of arrangements}

Geometric multidimensional continued fractions are defined for an arrangement of $n$ hyperplanes
passing through the origin (of dimension $n{-}1$) in $\r^n$ in general position.
Denote the set of all such arrangements by $\Omega_n$.

Any arrangement $\omega \in \Omega_n$ divides $\r^n$ into $2^{n}$ cones
(in fact, this is precisely the general position property for arrangements).
We say that these cones are the {\it cones of the arrangement $\omega$}.
Let us give now the following general definition.

\begin{definition}
A mapping $s$ of the coordinate cone $\r_{\ge 0}$ is called a {\it stretching}
of this cone if there exist positive real numbers $\alpha_1,\ldots, \alpha_n$
such that for every $(x_1,\ldots,x_n)\in \r_{\ge 0}$ we have
$$
s(x_1\ldots,x_n)=(\lambda_1 x_1,\ldots, \lambda_n x_n).
$$
Now let $C$ be a cone of some arrangement in $\Omega_n$.
We say that an automorphism $s$ of a cone $C$ is a {\it stretching} if
there exists an affine transformation taking $C$ to
$\r_{\ge 0}$ such that the the induced map is stretching for $\r_{\ge 0}$.

\vspace{1mm}

{\noindent
We say that two subsets $S_1$ and $S_2$ are {\it similar} if there exists a stretching of $C$ sending $S_1$ to $S_2$,
we write $S_1\sim S_2$ in this case.
}
\end{definition}

To be brief, we give the following general definition of geometric continued fractions
(that is a minor extension of a rather folklore definition which the author was communicated by A.~Ustinov).

\begin{definition}\label{mcf}
Let $C$ be a cone of some arrangement, and let $S$ be one of its subsets.
For an arbitrary subset $L\in C$ we define the set
$$
L_S=\bigcup \limits_{
U\cap L =\emptyset,
U\sim S
}
U.
$$

The boundary of $C\setminus L_S$ is called the {\it sail of $L\subset C$ with respect to $S$}.
We denote it by $\sail_{S,L}(C)$.
\end{definition}

There are several types of sets $S$ and $L$ traditionally considered in
various branches of geometry of numbers. Let us list them.

For $S$ it is natural to study the following sets:
\begin{itemize}
\item $\triangle^n$ is a basis $n$-simplex of the cone
(i.e., with a vertex at the origin, all other vertices on all distinct 1-edges of $C$);

\item $\square^n$ is a basis $n$-rectangular of the cone;

\item $\circ^n$ is an intersection of the unit circle centered at the vertex of $C$ (i.e. at the origin) with the cone $C$.

\item $h^n$ is a branch of a generalized hyperbola defined by
$$
\big|L_1(x_1,\ldots,x_n)\cdot \ldots \cdot L_n(x_1,\ldots,x_n)\big|=1,
$$
where $L_1,\ldots L_n$ are linear factors whose loci contain distinct hyperfaces of the cone $C$.
\end{itemize}

We continue with the following notation.
Let $\dot C$ be the punctured cone, i.e., the cone $C$ with its vertex being removed.

By $|\z^n|_C$ we denote the union of all the orbits of all integer points
with respect to the group of reflections about the hyperplanes of the arrangement containing $C$ that preserve all
the hyperplanes of the arrangement.

\vspace{2mm}

The following two sets are the natural choice for  $L$:
\begin{itemize}
\item $\z^n\cap \dot C$.

\item $|\z^n|_C\cap \dot C$.
\end{itemize}

\subsection{Important classes of geometric multidimensional continued fractions}

The sails for the above $S$ and $L$ deliver integer invariants of the arrangements in integer geometry.
Let us briefly describe the most famous multidimensional geometric generalizations of continued fractions.

\begin{example}
The classical two-dimensional case of geometric continued fractions is represented by sails
$$
\sail_{\triangle^2,\z^2\cap \dot C}(C).
$$
\end{example}

That is a particular case ($n=2$) of a general construction of multidimensional continued fractions
proposed by F.~Klein in 1895, see~\cite{Klein1895,Klein1896}.

\begin{definition}{\bf(Klein polyhedron.)}
{\it Klein polyhedron} of $C$ is the sail
$$
\sail_{\triangle,\z^n\cap \dot C}(C).
$$
\end{definition}
\vspace{2mm}

Let us compare this definition with the classical definition: {\it Klein polyhedron of a cone $C$
is the boundary of the convex hull of all integer points contained in $C$ except the vertex of the cone.}
It is clear that these two definitions agree for most interesting nondegenerate cases.
\vspace{2mm}

For further details regarding Klein continued fractions we refer
an interested reader to the books~\cite{Arnold2002} and~\cite{MCFbook}.

\begin{definition}{\bf(Minkowski-Voronoi polyhedron.)}
{\it Minkowski-Voronoi polyhedron} of $C$ is the sail
$$
\sail_{\square,|\z^n|_C\cap \dot C}(C).
$$
\end{definition}

This geometric generalization of ordinary continued fractions
was proposed independently by G.F.~Voronoi~\cite{Voronoi1896,Voronoi1952} and
H.~Minkowski~\cite{Minkowski1896,Minkowski1911}.
The algorithmic aspects of the construction of Minkowski-Voronoi continued fraction
can be found in~\cite{KU2017}.

\vspace{2mm}

{
\noindent
{\it Remark. }
{\bf$($Markov minima.$)$}
Theory of Markov minima of decomposable forms is also related to one of the cases of Definition~\ref{mcf}.
In fact the minima of such form are at vertices of sails
$$
\sail_{h^n,\z^n\cap \dot C}(C_i), \quad i=1,\ldots, 2^n,
$$
where $C_i$ are all the cones of the arrangement of hyperspaces that are the kernels for the linear factors of the form.

For the classical case of $n=2$ we refer a reader to
the original papers of Markov~\cite{Markoff1879, Markoff1880}, and to a nice overview about Markov spectra in~\cite{Cusic1989} by T.~Cusick and M.~Flahive.
}

\vspace{2mm}

{
\noindent
{\it Remark. }
{\bf$($Hermite's parametrisation.$)$}
The sails
$$
\sail_{\circ^n,\z^n\cap \dot C}(C)
$$
are used in the Method of Hermit's parametrisation~\cite{Ozhigova1982}.
}

\section{Open problems on geometric multidimensional continued fractions}

One of the main source of problems in this area are the problems of face structure description of sails and the corresponding statistics. Indeed sails are either polyhedral surfaces
(in cases of $\triangle^n$ and $\square^n$) or a little more complicated stratified surfaces with singularities (in cases of $\circ^n$ and $h^n$).
In both cases we have a combinatorial structure with some lattice invariants involved.

\subsection{Directions for further study in general}

There are four essentially different cases of geometric multidimensional continued fractions to study:

\vspace{2mm}

\begin{itemize}
\item[{[}R{]}] {\bf Rational cones:} These are the cones that contain integer points at each of the 1-edge of the cone. The sails of such cones contain finite number of faces of all dimensions.

\vspace{2mm}

\item[{[}A{]}] {\bf Algebraic cones:} Consider a matrix $A\in \GL_n(\z)$, i.e., all the elements of $A$ are integers and the determinant of $A$ is $\pm 1$.
In addition we assume that characteristic polynomial of $A$ is irreducible over $\q$ and has only real roots ({\it totally real case}).
Then the arrangement defined by the invariant planes of $A$ is called {\it algebraic}.
It's sails are {\it algebraic} as well. By Dirichlet unit theorem, there is a transitive action of the Dirichlet subgroup for $A$.
Recall that the Dirichlet subgroup for $A$ is the group of all $\SL_n(\z)$-matrices commuting with $A$.
In the totally real case, the Dirichlet subgroup for $A$ is isomorphic to $\z^{n-1}$.
And therefore, the factor-sails with respect to Dirichlet subgroup of $A$ are homeomorphic to the $(n{-}1)$-dimensional torus.
For each sail, the corresponding torus is subdivided into faces equipped with lattice invariants of the sail.

\vspace{2mm}

\item[{[}D{]}] {\bf Non-periodic degenerate cones:} Here are the cones that contain integer points on their faces.
Usually the sail of such cones are neither periodic nor finite.
Furthermore, as a rule such sail are not combinatorial polyhedra.
The sails for non-periodic degenerate cones serve as counterexamples for several questions.

\vspace{2mm}

\item[{[}G{]}] {\bf Non-periodic generic cones:} These are non-algebraic cones with no integer points on faces. Such cones usually possess non-periodic sails that are combinatorial polyhedra (see~\cite{MCFbook} and~\cite{Lachaud2002} for more details).
\end{itemize}

\vspace{2mm}

The following two question form the core of combinatorial and statistical problems for different sails.
\begin{itemize}
\item Find an adequate complete invariant of polyhedra that appear as faces at sails.

\vspace{1mm}

\item Find an adequate complete invariant of cones up to integer congruences.

\vspace{1mm}

\item Study the face distribution of sails (generalized Gauss-Kuzmin distribution).
\end{itemize}

In next three subsection we briefly address these major directions.

\subsection{Faces of sails}

The central question of this subsection is as follows.
\begin{problem}
Classify all combinatorial possible types of faces.
\end{problem}
This question targets mostly Klein polyhedra
(i.e., $\sail_{\triangle^n,\z^n\cap \dot C}(C)$).
In relation to this question various view-obstacle problems arise.
One of the first open famous problems here is generalization of White's theorem,
which we discuss below.

\vspace{2mm}

A polyhedron $P$ is said to be {\it empty} if its intersection with $\z^n$ is precisely the set of vertices of $P$.
We have the following problem.

\begin{problem}
Classify all empty simplices of dimension $n$ up to lattice congruence.
\end{problem}

If $n=1,2$ then there exists a unique empty simplex up to integer congruences,
it is a coordinate simplex.
For the case $n=3$ there is a family of empty simplices,
they are classified in White's theorem (see~\cite{White1964}).
The cases of $n>3$ remain open.

\vspace{2mm}

The next problem here is as follows.
\begin{problem}
Which $n$-gons are realizable as faces of an $m$-dimensional continued fraction?
Here are two essentially geometrically different subcases:

\begin{itemize}
\item {\bf Faces at integer distance 1 to the origin:}
this is a question of description of integer convex polyhedra that are inscribed to simplices of full dimension in $\r^n$.

\item {\bf Faces at integer distance greater than 1 to the origin:} here we have a typical view-obstacle problem, one avoids integer points that are view-obstacles for the corresponding faces.
\end{itemize}
\end{problem}

Both statements of the theorem are trivial for $n=2$. 
Two-dimensional faces on integer distance greater than $1$ 
are classified in~\cite{Karpenkov2005}).
All the other cases remain open.

\vspace{2mm}

Finally in the context of this subsection we would like to mention a classical Lonely Runner conjecture 
(for more details and discussions we refer to~\cite{MCFbook}).

\begin{problem}{\bf( Lonely Runner conjecture.)} Assume that
$k$ runners start a competition at the a common point of the circular track with circumference 1. Let them all have distinct constant speeds.
Then for any given runner, there is a time at which that runner is lonely, 
i.e., at arc distance at least $1/k$ away from every other runners.
\end{problem}

\subsection{Combinatorial structure of sails}

One can easily define tautological complete invariant 
which does not show the combinatorics of sails. 
Such invariants are not informative.
Most of the known descriptive invariants are
\\
--- either {\bf not complete}, i.e., they do not distinguish some of integer non-congruent objects;
\\
--- or {\bf not precise}, i.e., the set of configurations of this invariant that can be realized as faces or cones is not known.

\vspace{2mm}

So the following problem remains open for Klein sails, Minkowski-Voronoy polyhedra, and most other continued fractions (for the case [R] as above).
\begin{problem}
Describe all finite two-dimensional sails (and the corresponding continued fractions).
\end{problem}

The following problem is on finite continued fractions contained in the different angles of a tetrahedron.
That problem is a version of a question on singularities of toric varieties in algebraic geometry and on cuspidal singularities in singularity theory
(compare with Problem~\ref{IKEA2d} above).
\begin{problem}{\bf (Multidimensional IKEA problem.)}
Describe the collections of the sails of the cones for all polytopes of a given combinatorial type.
\end{problem}

Most of further questions are related to the algebraically periodic
multidimensional sails and their fundamental domains (case [A]).
We begin with the following question is of interest.
\begin{problem}{\bf (V.~Arnold.)}
Does there exist an algorithm to decide whether a given type of fundamental domain is realizable by a periodic continued fraction?
\end{problem}
This problem is relevant for various types of sails (including Minkowski-Voronoi and Klein cases).

\vspace{2mm}

{
\noindent
{\it Remark.}
For the case of Klein sails (i.e., $\sail_{\triangle,\z^n\cap \dot C}(C).$) 
it is still unknown if there exists a two-dimensional periodic sail which contains only quadrangles (even a single quadrangle in the fundamental domain).
}

\vspace{2mm}

A {\it torus decomposition corresponding to a periodic sail in the sense of Klein} is the collections of integer congruence types of polyhedra in any  fundamental domain endowed with the integer distances to the vertex of the cone.
\begin{conjecture}{\bf (V.~Arnold.)}
Torus decompositions of integer noncongruent Klein sails are distinct.
\end{conjecture}

Next open problem is considered as a complete generalization of the Lagrange theorem
on periodicity of ordinary continued fractions.
\begin{problem} {\bf (V.~Arnold.)}
Describe all torus decompositions that are realized by periodic two-dimensional continued fractions.
\end{problem}

The following three observations are confirmed by numerous examples.
\begin{conjecture}
The torus decomposition of every algebraic sail in the sense of Klein
contains at least one face that is
\begin{itemize}
\item integer congruent to a triangle;

\item at the unit integer distance form the vertex of the cone;

\item at the integer distance from the vertex of the cone greater than 1.
\end{itemize}
\end{conjecture}


In the theory of field extensions the following problem is actual.
\begin{problem} {\bf (V.~Arnold.)}
Classify continued fractions that correspond to the same cubic extension of the field of rational numbers.
\end{problem}
In fact, this problem seems to be open even for the classical case of ordinary periodic continued fractions.

\vspace{2mm}

We would like to conclude this subsection with a very general problem of construction generic cones.
\begin{problem}
Prove the existence of a cone for a single non-periodic 
combinatorial structure ($n\ge 3$).
\end{problem}

This question is really trivial for standard continued fractions $(n=2)$. 
For instance there exists a cone, whose continued fraction is
$$
[1;2:3:4:\ldots].
$$
Currently there is no examples of non-periodic multidimensional non-periodic sails
with a prescribed structure. 
This concerns the case [G] for both $\z^n\cap \dot C$
and $|\z^n|\cap \dot C$ and for all choices of $S$ (including $\triangle^n$, $\circ^n$, $\square^n$, $h^n$). The choice of $h^n$ is somehow different from the rest geometrically, nevertheless there is no results known to the author regarding the last problem.

\subsection{Sail statistics}

Recall that Gauss-Kuzmin statistics in the two dimensional case is as follows. 
For every $k>0$ and $n>0$ the probability of an event $a_n=k$
in a generic continued fractions is given by the following formula:
$$
\frac{1}{\ln 2}\ln\Big(1-\frac{1}{(k+1)^2}\Big).
$$
Surprisingly, the last expression has a remarkable interpretation in terms of projective geometry
$$
\frac{1}{\ln 2}\ln\Big(1-\frac{1}{(k+1)^2}\Big)=
\frac{1}{\ln 2}\ln\Big({k(k+2)}{(k+1)^2}\Big)=
\frac{\ln\big([-1,0,k,k{+}1]\big)}{\ln\big([-1,0,1,\infty]\big)},
$$
where $[a,b,c,d]$ is the cross-ratio for $a$, $b$, $c$, and $d$.
This is not a coincidence as geometric continued fractions are invariant with respect to projective group acting on two-dimensional geometric continued fractions which roughly coincide with ordinary continued fractions (see~\cite{Karpenkov2007} for more details).

Note that the cross-ratio in the denominator is the sum of probabilities 
of all element ``$k$'' for $k=1,2,3,\ldots$, so we have
$$
\ln\big([-1,0,1,\infty]\big)=\sum\limits_{k=1}^{\infty}\ln\big([-1,0,k,k{+}1]\big).
$$

\vspace{2mm}

On the one hand there is not much known for the multidimensional generalization of 
Gauss-Kuzmine statistics, as there is no straightforward generalization 
of Gauss map that is related to multidimensional sails. 
On the other hand there is a geometric approach to define geometric relative frequencies by means of hyperbolic geometry. Here one integrates infinitesimal cross-ratio with respect to the unique M\"obius measure on the space of all arrangements, which is invariant under projective transformations acting on it.
We refer an interested reader to~\cite{Illarionov2012,Karpenkov2007,MCFbook}.

Here we would like to outline the following three questions (which are open for $n\ge 3$).
\begin{problem}
Find frequencies on $n$-dimensional continued fractions with the highest relative frequencies.
\end{problem}

\begin{problem}
For every positive integer constant $C$ there exist only finitely many pairwise integer non-congruent faces with frequencies exceeding $C$.
\end{problem}

\begin{problem}
Is that true that sum of all relative frequencies for all possible faces is finite for higher dimensions $(n\ge 3)$?
\end{problem}
\begin{problem}
In case of positive answer to the above question find the generalization of the Gauss map and compare the corresponding frequencies of faces with the related frequencies coming from M\"obius geometry.
\end{problem}

\section{Further open questions}

In this section we briefly mention some directions that concern different geometric approaches to multidimensional continued fractions. For further details we refer to~\cite{MCFbook}.

The first question is to generalize Farey summation and its representation as a Farey tessellation.
\begin{problem}
Find a natural generalization of the Farey tessellation to higher-dimensional hyperbolic geometry.
\end{problem}

The second question concerns periodicity of Jacobi-Perron algorithms.
\begin{problem}{\bf (Jacobi's last theorem.)}
Let $K$ be a totally real cubic number field.
Consider arbitrary elements $y$ and $z$ of $K$ such that $0<y,z<1$ (here we assume that $1$, $y$, and $z$ are independent over $\mathbb \q$).
Is it true that the Jacobi-Perron algorithm generates an eventually periodic continued fraction with starting data $v=(1,y,z)$.
\end{problem}
An essential contribution to this problem was made in~\cite{Assaf2005}.

\vspace{2mm}

The third question is related to the study of minima of decomposable forms.
\begin{problem}
Study geometric properties of Markov spectrum.
\end{problem}

\vspace{2mm}

Finally, we would like to address the following question on low-dimensional topology.
Recall that ordinary finite continued fractions provide a complete invariant to rational knots, known also as two-bridge knots.
So the following rather classical question is actual.
\begin{problem}
Generalize continued fractions to describe 3-bridge knots.
\end{problem}

\vspace{2mm}

{\noindent
{\bf Acknowledgement.}
The author is grateful to A.~Ustinov for communicating a general definition of geometric multidimensional continued fractions.
}


\vspace{5mm}

\end{document}